\documentclass{article}

\usepackage{arxiv}

\usepackage{amsthm,amsmath,amssymb,a4wide,xcolor,graphicx,tikz, hyperref,physics}

\usepackage[utf8]{inputenc} 
\usepackage[T1]{fontenc}    
\usepackage{hyperref}       
\usepackage{url}            
\usepackage{booktabs}       
\usepackage{amsfonts}       
\usepackage{nicefrac}       
\usepackage{microtype}      
\usepackage{cleveref}       
\usepackage{lipsum}         
\usepackage{graphicx}
\usepackage{doi}

\usepackage[fontsize=11pt]{fontsize}

\newcommand\slo{s}
\newcommand\nlc{\eta}

\newcommand\phih{\hat{\phi}}
\newcommand\psih{\hat{\psi}}

\newcommand\zbot{y}
\definecolor{grey}{rgb}{0.5, 0.5, 0.5}

\newcommand{\uproman}[1]{\uppercase\expandafter{\romannumeral#1}}
\newcommand{\lowroman}[1]{\romannumeral#1\relax}
\title{A paraxial approach for the inverse problem of vibro-acoustic imaging in frequency domain}

\newif\ifuniqueAffiliation
\uniqueAffiliationtrue

\ifuniqueAffiliation 
\author{ 
	\hspace{1mm}Teresa Rauscher \\
	Department of Mathematics\\
	University of Klagenfurt\\
	Klagenfurt, Austria \\
	\texttt{teresa.rauscher@aau.at} \\
}
\else
\usepackage{authblk}

\affil[1]{Department of Computer Science, Cranberry-Lemon University, Pittsburgh, PA 15213}
\affil[2]{Department of Electrical Engineering, Mount-Sheikh University, Santa Narimana, Levand}
\fi

\renewcommand{\headeright}{}
\renewcommand{\undertitle}{}
\renewcommand{\shorttitle}{Paraxial vibro-acoustic imaging}

\hypersetup{
pdftitle={A template for the arxiv style},
pdfsubject={q-bio.NC, q-bio.QM},
pdfauthor={David S.~Hippocampus, Elias D.~Striatum},
pdfkeywords={First keyword, Second keyword, More},
}

\begin{document}
\maketitle

\begin{abstract}
	This paper presents a paraxial modeling approach for vibro-acoustography, a high-frequency ultrasound imaging technique that makes use of the excited low-frequency field to achieve a higher resolution while avoiding speckles. We start from a general second order wave equation, introduce a second order perturbation and make use of a paraxial change of variables to pass over to directed beams. Depending on the relation between the second order perturbation and the directivity of the beams this leads to different systems of PDEs that involve a spatially variable parameter governing the interaction of the incoming beams to generate the propagating acoustic field. We then consider the inverse problem of determining this parameter for one of these models where we present a Landweber-Kaczmarz reconstruction method to obtain a spatial image of the region of interest.
\end{abstract}


\section{Introduction}
In medical diagnostics the quest for non-invasive, safe and precise imaging techniques has been relentless. While traditional ultrasound imaging is well-established and widely used, there is a continuing interest in nonlinear acoustics as the demand for higher image resolution and diagnostic accuracy continues to grow. However, nonlinear models also come with drawbacks such as scattering from small inclusions and stronger attenuation at higher frequencies. In this paper we consider vibro-acoustic imaging by means of ultrasound imaging \cite{FatemiGreenleaf:1998,FatemiGreenleaf:1999}. This technique avoids the drawbacks of nonlinear models while achieving speckle-free images of high resolution and providing valuable diagnostic information on tissue properties.  

In contrast to conventional ultrasound imaging, where a beam of ultrasound waves is sent into the tissue using a transducer, in vibro-acoustic imaging two ultrasound beams at high but slightly different frequencies are sent into the medium.  At a focus they interact nonlinearly and this interaction excites a wave field that basically propagates at the difference frequency and is eventually measured by a hydrophone. For more details on the experiments see \cite{FatemiGreenleaf:1998,FatemiGreenleaf:1999} and for a graphical illustration of the described experiment for image acquisition see Figure \ref{fig: exp}. 

Vibro-acoustography can be used in many biomedical and non-destructive applications, especially, to image several types of tissues, such as breast, liver, prostate, arteries, and thyroid. Typically, the ultrasound frequencies used are in the range of $1$ to $10$ MHz and the difference frequencies range form $10$ to $100$ kHz. For wide applicability in clinical settings one uses linear-array transducers for forming vibro-acoustography beams. More information on beamforming and implementation on clinical ultrasound systems can be found in \cite{UrbanFatemi:2013, UrbanFatemi:2011, Multifrequ:2006}. 

A modeling and simulation framework for vibro-acoustic imaging has been presented in \cite{MalcolmReitichYangGreenleafFatemi:2008:parabolic, MalcolmReitichYangGreenleafFatemi:2008:Numerics}. The model consists of a system of partial differential equations that involve the speed of sound $c=c(x)$ and the nonlinearity parameter $\frac{B}{A}=\frac{B}{A}(x)$ which depend on space due to the inhomogeneity of the medium. Therefore, their reconstruction yields a spatial image of the region of interest. In case of reconstructing the individual coefficients $c=c(x)$ or $\frac{B}{A}=\frac{B}{A}(x)$, this is related to ultrasound tomography and nonlinearity parameter imaging, respectively, cf. e.g. \cite{Bjorno:1986,Cain:1986,IchidaSatoLinzer:1983} and the citing literature. The corresponding inverse problem for vibro-acoustic imaging has been put into a mathematical framework in \cite{Kaltenbacher:2022}. In this model all waves are assumed to propagate omni-directional, although the high frequency waves actually have a strongly preferred direction of propagation, which can justify the use of a parabolic approximation, cf. e.g. \cite{Tappert:1977}. 

There are several methods of deriving parabolic approximations \cite{Jensen:2011}. In \cite{MalcolmReitichYangGreenleafFatemi:2008:parabolic, MalcolmReitichYangGreenleafFatemi:2008:Numerics} a diagonalization procedure has been chosen for the parabolic approach and for efficient numerical simulation a decomposition approach has been devised. The forward problem is split into three components: (\lowroman{1}) directed high frequency propagation of the two ingoing beams $\phi_1$, $\phi_2$, (\lowroman{2}) nonlinear interaction of these at the focus, and (\lowroman{3}) undirected low frequency propagation of $\psi$ to the hydrophone. 

In this paper the goal is to derive directed models making use of the paraxial approximation to derive the Khokhlov-Zabolotskaya-Kuznetsov (KZK) equation \cite{Zabolotskaya:1969} from Kuznetsov's equation. The KZK equation is a directed nonlinear wave equation that is particularly important for modeling and understanding the behavior of ultrasound waves in medical ultrasound applications as it takes various effects including nonlinearity, diffraction and absorption into account, see also, e.g. \cite{Rozanova:2007} for the analysis of the KZK equation. As we pass from undirected beams to directed beams, the equations change significantly and an additional small parameter for the order of the direction comes into play.  
 
Throughout this paper we assume the speed of sound to be constant and only the nonlinearity parameter to be space dependent. Reconstructing both parameters simultaneously the paraxial change of variables has to be adapted and redefined, see equation \eqref{parax_deriv}, since this transformation of variables changes the time variable to the retarded time, which means that the time gets space dependent and therefore, also dependent on the speed of sound. Results on the simultaneous reconstruction of sound speed and nonlinearity parameter can be found in \cite{2023simultaneous}.

The plan of this paper is as follows. In Section \ref{sec: model}, we first of all derive several paraxial models starting from the second order wave equation, introducing a second order perturbation for the velocity potential and then performing the paraxial change of variables. Then we will transfer the time domain formulations to frequency domain, equip the system with suitable boundary conditions and give some remarks on adaptions and generalizations of the model and its experimental setup. Focusing on one of these models in Section \ref{sec: IP}, we will state the inverse problem of nonlinearity imaging and present a Landweber-Kaczmarz method to solve the linear inverse problem of reconstructing the nonlinearity parameter. Algorithms for the direct problem and also the inverse problem can be found in Section \ref{sec: num}.


\newpage

\section{The model} \label{sec: model}

\begin{figure}[h]
	\centering
	\includegraphics[scale=0.46]{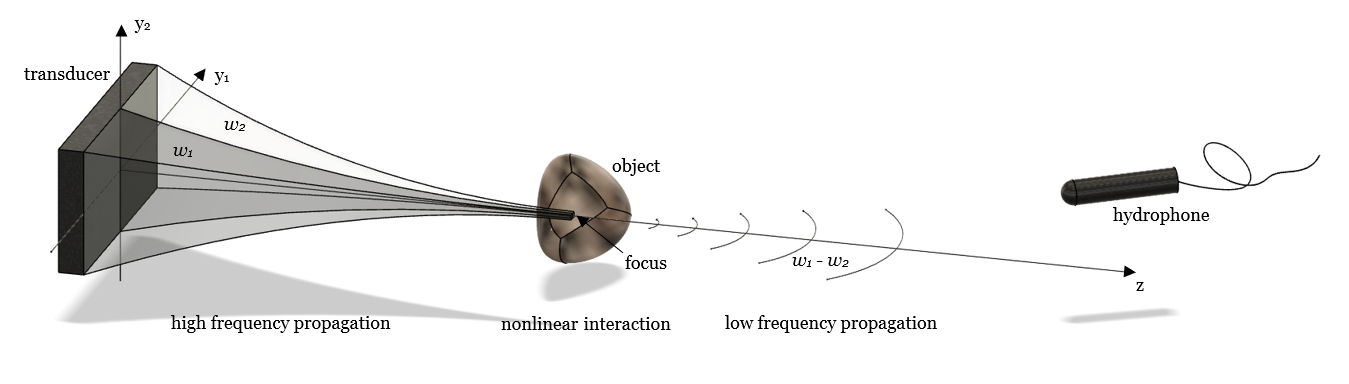} 
	\caption{experimental setup}
	\label{fig: exp}
\end{figure}

In vibro-acoustography two ultrasound beams with velocity potentials $\phi_1$, $\phi_2$ are excited by a (linear) transducer. They are focused at two slightly different high frequencies $\omega_1$, $\omega_2$ ($\Delta \omega = \omega_1-\omega_2 \ll \omega_1,\, \omega_2$) with $\omega_1 > \omega_2 > 0$. These ingoing beams have a strong preferred direction of propagation, here on the central axis $z$, away from the transducer. At a focus they interact nonlinearly and this interaction excites a wave field with velocity potential $\psi$ that basically propagates at the difference frequency $\Delta \omega$, which is then quite low. Since the source of the low frequency field is point-like the resulting field is omni-directional. The acoustic response is eventually measured by a hydrophone located far from the focus point and an image can be formed. The experimental setup for image acquisition is illustrated in Figure \ref{fig: exp}.

We assume that the beam is formed by a linear transducer. Note that there are different ways to arrange the two beams. We choose to split the transducer array in left and right, since we want to keep the beams symmetrically to the central axis and the focus. This is emphasized in Figure \ref{fig: exp} by the frequency $\omega_1$ on the left and $\omega_2$ on the right. Some further remarks on adaptions and generalizations of the experimental setup can be found in Section \ref{sec: adaptions}.

\subsection{Derivation of paraxial models}

We now want to look at the derivation of a directive model for the ingoing beams in this methodology. We assume that the beams are naturally focused, so by geometrical means and make use of the parabolic approximation for the derivation of the KZK equation from Kuznetsov's equation, cf. e.g., \cite{Rozanova:2007}. 

\begin{figure}[h]
	\centering
	\begin{minipage}{\textwidth}
		\hspace{1.5cm}
		\begin{minipage}{0.45\textwidth}
			\begin{tikzpicture}
				\coordinate (O) at (0,0,0);
				\draw[thick, ->] (O) -- (3,0,0) node[anchor=north east]{$x_1$};
				\draw[thick, ->] (O) -- (0,0,3) node[anchor=north east]{$x'$};
				\draw[thick, ->] (O) -- (0,2.3,0) node[anchor=north west]{$t$};
			\end{tikzpicture}
		\end{minipage} 
		\hspace{-1cm}
		$\Rightarrow$
		\hspace{0.7cm}
		\begin{minipage}{0.45\textwidth}
			\begin{tikzpicture}
				\coordinate (O) at (0,0,0);
				\draw[thick, ->] (O) -- (3,0,0) node[anchor=north]{$z=\tilde{\varepsilon} x_1$};
				\draw[thick, ->] (O) -- (0,0,3) node[anchor=north]{$y=\sqrt{\tilde{\varepsilon}}x'$};
				\draw[thick, ->] (O) -- (0,2.3,0) node[anchor=north west]{$\tau=t-\frac{x_1}{c}$};
			\end{tikzpicture}
		\end{minipage}
	\end{minipage}
\caption{paraxial change of variables}
\label{fig: paraxial}
\end{figure}
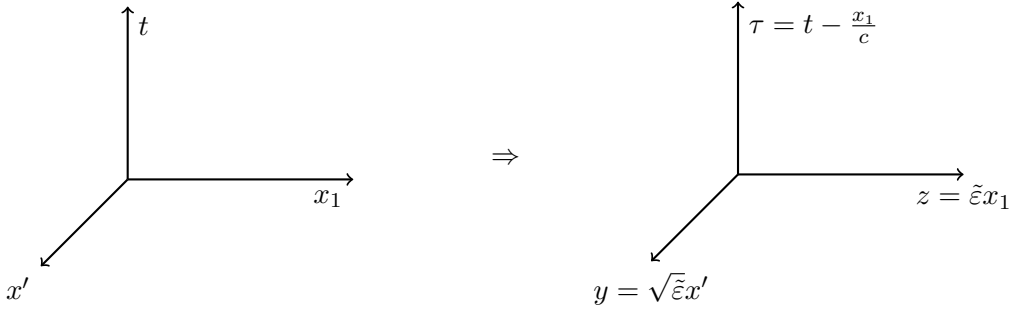

We assume that the direction of propagation is the $x_1$-axis. Then, the acoustical properties of directive beam's propagation are that the beams
\begin{itemize}
	\item are concentrated near the $x_1$-direction;
	\item propagate along the $x_1$-direction;
	\item are generated either by an initial condition or by a forcing term on the boundary $x_1=0$.
\end{itemize}
 Furthermore, it is assumed that the variation of beam's propagation in the direction $x'=(x_2, x_3, \ldots, x_n)$ perpendicular to the $x_1$-axis is much larger than its variation along the $x_1$-axis. This can mathematically be modeled by a paraxial change of variables. For the initial variables $(t, x_1, x')$, where $x'=(x_2, \ldots, x_n)$, we perform the paraxial change of variables 
\begin{align} \label{parax}
	(\tau, z,y)=\left( t-\frac{x_1}{c}, \tilde{\varepsilon} x_1, \, \sqrt{\tilde{\varepsilon}} x'\right),
\end{align} 
where $\tilde{\varepsilon}$ is a small parameter. This change of variables is illustrated in Figure \ref{fig: paraxial}. \\

\vspace{0.2cm}
In the paraxial change of variables in equation \eqref{parax} the speed of sound is assumed to be constant. If considered it variable the definition \eqref{parax}, in particular the retarded time variable, would have to be adapted such that
\begin{align}\label{parax_c}
	(\tau, z, y)=\left( t-\check{\mathfrak{t}}, \tilde{\varepsilon} x_1, \sqrt{\tilde{\varepsilon}}x'\right), 
\end{align}
with $\check{\mathfrak{t}}=\check{\mathfrak{t}}(x_1,x'), \, \mathfrak{t}=\mathfrak{t}(z,\zbot)$, where 
\begin{equation*}
	\frac{1}{c^2}=|\nabla_x\check{\mathfrak{t}}|^2
	=|\partial_{x_1}\check{\mathfrak{t}}|^2+|\nabla_{x'}\check{\mathfrak{t}}|^2
	=\tilde{\varepsilon}^2|\partial_{z}\mathfrak{t}|^2+\tilde{\varepsilon}|\nabla_{\zbot}\mathfrak{t}|^2, 
	\quad \check{\mathfrak{t}}(0,x')=0, \ x'\in \tilde{\Omega}_{\zbot}.
\end{equation*}

In the following, we assume constant speed of sound and consider the definition for the paraxial change of variables in \eqref{parax}. 

To derive a model for vibro-acoustic imaging, we start with the general wave equation in terms of the velocity potential $\phi$ given by
\begin{equation} \label{eq: der15}
	\partial_t^2 \phi - c^2 \nabla^2 \phi = \partial_t |\nabla \phi |^2 - \frac{1}{2} (\nabla \phi) \cdot \nabla |\nabla \phi|^2 + (\gamma -1) \partial_t \phi \nabla^2 \phi - \frac{\gamma -1}{2} |\nabla \phi |^2 \nabla^2 \phi,
\end{equation}
where $\gamma=\frac{B}{A}+1$ and $\phi=\phi(x,t)$ with $x=(x_1,\ldots, x_n) \in \mathbb{R}^n$, $t \in \left[ 0,T\right]$ for $T>0$. Equation \eqref{eq: der15} is exact for a perfect lossless gas and often serves as a starting point for perturbation analyses in nonlinear acoustics. Approximations lead to simplified wave equations, e.g. the Westervelt equation, that are more amenable to analysis. For more details and the derivation see \cite{Hamilton:1998}. 

Since the acoustic state variables are all $O(\varepsilon)$ or smaller, one can expect that the velocity potential depends analytically on $\varepsilon$. Therefore, we introduce the following second order approximation for the velocity potential $\phi$ such that
\begin{equation} \label{eq: dervib0}
	\phi(x,t,\varepsilon) = \varepsilon(\phi_1(x,t)+\phi_2(x,t)) + \varepsilon^2 \psi(x,t),
\end{equation}
where $\varepsilon$ is very small, such that $\phi_1$, $\phi_2$ represent the high-frequency beams and $\psi$ the low-frequency field. 

Now, plugging this power series expansion of $\phi$ relative to $\varepsilon$ into the general wave equation given by equation \eqref{eq: der15} yields 
\begin{align*} 
	\varepsilon\partial_t^2 (\phi_1+\phi_2) + \varepsilon^2 \partial_t^2 \psi &-c^2 \varepsilon \nabla^2(\phi_1+\phi_2)-c^2 \varepsilon^2 \nabla^2 \psi \\
	& = \varepsilon^2\partial_t |\nabla(\phi_1+\phi_2)|^2 + \varepsilon^2(\gamma-1) \partial_t (\phi_1+\phi_2) \nabla^2(\phi_1+\phi_2) + O(\varepsilon^3).
\end{align*} \normalsize
Considering only terms up to order $\varepsilon^2$ and rearranging the equation above gives 
\begin{align} \label{eq: dervib1}
	\varepsilon \left[ \partial_t^2 (\phi_1+\phi_2) -c^2 \nabla^2(\phi_1+\phi_2) - \varepsilon \partial_t |\nabla(\phi_1+\phi_2)|^2 -  \right. & \left. \varepsilon (\gamma-1) \partial_t (\phi_1+\phi_2) \nabla^2(\phi_1+\phi_2) \right]  \notag \\
	& = \varepsilon^2\left[ - \partial_t^2 \psi + c^2  \nabla^2 \psi\right] .
\end{align} \normalsize

One can observe that we arrived at well known wave equations. The equation on the left hand side in terms of the velocity potential $\phi_1+\phi_2$ is known as the Kuznetsov equation (without viscosity) and the right hand side is just the linear wave equation in terms of $\psi$. 

The next step is to include the concept of directive beams in this model. Therefore, we perform the paraxial change of variables \eqref{parax}, 
where $\tilde{\varepsilon} \ll 1$, so $\tilde{\varepsilon} \ll \sqrt{\tilde{\varepsilon}}$ and in general $\tilde{\varepsilon}$ can be different from $\varepsilon$ that was introduced for the second order approximation in \eqref{eq: dervib0} but for the model to make sense we will see that the relation $\varepsilon \leq \tilde{\varepsilon} \leq \sqrt{\varepsilon}$ should hold. For the derivatives this transform of variables yields
\begin{equation}\label{parax_deriv}
	\partial_t=\partial_{\tau}, \quad \nabla=\left( -\frac{1}{c} \partial_{\tau} + \tilde{\varepsilon} \partial_z, \sqrt{\tilde{\varepsilon}} \nabla_y\right) , \quad \nabla^2=\frac{1}{c^2} \partial_{\tau}^2 - 2 \frac{\tilde{\varepsilon}}{c} \partial^2_{\tau z} + \tilde{\varepsilon}^2 \partial_z^2 + \tilde{\varepsilon} \nabla_y^2
\end{equation}
and for the term $|\nabla \phi|^2$ one arrives at
\begin{equation*}
	|\nabla \phi|^2 = \nabla \phi \cdot \nabla \phi = \frac{1}{c^2}(\partial_{\tau} \phi)^2 - 2 \frac{\tilde{\varepsilon}}{c} \partial_{\tau} \phi \partial_z \phi + \tilde{\varepsilon}^2 (\partial_z \phi )^2 + \tilde{\varepsilon} |\nabla_y \phi|^2.
\end{equation*} 
For the left hand side of equation \eqref{eq: dervib1} this change of variables then leads to 
\begin{align*}
	\varepsilon &\left[ \partial_t^2 (\phi_1+\phi_2) -c^2 \nabla^2(\phi_1+\phi_2) - \varepsilon \partial_t |\nabla(\phi_1+\phi_2)|^2 -   \varepsilon (\gamma-1) \partial_t (\phi_1+\phi_2) \nabla^2(\phi_1+\phi_2) \right]  \notag \\
	&= \varepsilon  \left[ 2 \tilde{\varepsilon}c \partial^2_{\tau z} (\phi_1+\phi_2)- \tilde{\varepsilon}^2 c^2 \partial_z^2 (\phi_1+\phi_2)- \tilde{\varepsilon} c^2\nabla_y^2 (\phi_1+\phi_2) \right.  \\
	& \left. - \varepsilon \left( \frac{\gamma+1}{2c^2} \partial_{\tau}\left| \partial_{\tau} (\phi_1+\phi_2)\right|^2 
	- 2 \frac{\tilde{\varepsilon}}{c} \left( \partial_{\tau}^2 (\phi_1+\phi_2)\partial_{z} (\phi_1+\phi_2)+ \gamma \partial_{\tau}(\phi_1+\phi_2) \partial_{\tau z}^2 (\phi_1+\phi_2) \right) \right. \right. \\
	& \left. \left.  + \tilde{\varepsilon}^2 \left( \partial_{\tau}  (\partial_z (\phi_1+\phi_2) )^2 + (\gamma-1) \partial_{\tau} (\phi_1 + \phi_2) \partial_z^2(\phi_1+\phi_2) \right) \right. \right. \\
	& \left. \left. + \tilde{\varepsilon} \left(  \partial_{\tau}  |\nabla_y (\phi_1+\phi_2)|^2 +(\gamma-1) \partial_{\tau}(\phi_1+\phi_2) \nabla_y^2 (\phi_1+\phi_2)\right) \right)\right].
\end{align*} \normalsize

A paraxial change of variables for the velocity potential $\psi$ for the right hand side of equation \eqref{eq: dervib1} yields
\begin{align*}
	\varepsilon^2\left[ -  \partial_t^2 \psi + c^2 \nabla^2 \psi\right]  
	& = \varepsilon^2 \left[ -2 \tilde{\varepsilon} c \partial_{\tau z} \psi + c^2 \tilde{\varepsilon}^2 \partial_z^2 \psi + \tilde{\varepsilon} c^2 \nabla_y^2 \psi \right].
\end{align*} 
\normalsize

Therefore, altogether we end up with the equation 
\begin{align*}
	\varepsilon & \left[ 2 \tilde{\varepsilon}c \partial^2_{\tau z} (\phi_1+\phi_2)- \tilde{\varepsilon}^2 c^2 \partial_z^2 (\phi_1+\phi_2)- \tilde{\varepsilon} c^2\nabla_y^2 (\phi_1+\phi_2) \right.  \\
	& \left. - \varepsilon \left( \frac{\gamma+1}{2c^2} \partial_{\tau}\left| \partial_{\tau} (\phi_1+\phi_2)\right|^2- 2 \frac{\tilde{\varepsilon}}{c} \left( \partial_{\tau}^2 (\phi_1+\phi_2)\partial_{z} (\phi_1+\phi_2)+ \gamma \partial_{\tau}(\phi_1+\phi_2) \partial_{\tau z}^2 (\phi_1+\phi_2) \right) \right. \right. \\
	& \left. \left.  + \tilde{\varepsilon}^2 \left( \partial_{\tau}  (\partial_z (\phi_1+\phi_2) )^2 + (\gamma-1) \partial_{\tau} (\phi_1 + \phi_2) \partial_z^2(\phi_1+\phi_2) \right) \right. \right. \\
	& \left. \left. + \tilde{\varepsilon} \left(  \partial_{\tau}  |\nabla_y (\phi_1+\phi_2)|^2 +(\gamma-1) \partial_{\tau}(\phi_1+\phi_2) \nabla_y^2 (\phi_1+\phi_2)\right)    \right)\right]  \\
	= & \, \varepsilon^2 \left[ -2 \tilde{\varepsilon} c \partial_{\tau z} \psi + c^2 \tilde{\varepsilon}^2 \partial_z^2 \psi + \tilde{\varepsilon} c^2 \nabla_y^2 \psi \right].
\end{align*}

Reordering with respect to the different orders of $\varepsilon$ and $\tilde{\varepsilon}$ yields 
\begin{align} 
	\varepsilon \tilde{\varepsilon} \left[  2 c \partial^2_{\tau z} (\phi_1+\phi_2) -c^2 \nabla_y^2 (\phi_1+\phi_2) \right]  &=0, \label{eq: sys1} \\ 
	\varepsilon^2 \left[  -\frac{\gamma +1}{2c^2} \partial_{\tau}\left| \partial_{\tau} (\phi_1+\phi_2)\right|^2 \right]   &=0, \label{eq: sys2}\\ 
	\varepsilon \tilde{\varepsilon}^2 \left[  -c^2 \partial_z^2 (\phi_1+\phi_2) \right]  & =0, \label{eq: sys3} \\
	\varepsilon^2 \tilde{\varepsilon} \left[ \frac{2}{c} \left( \partial_{\tau}^2 (\phi_1+\phi_2)\partial_{z} (\phi_1+\phi_2)+ \gamma \partial_{\tau}(\phi_1+\phi_2) \partial_{\tau z}^2 (\phi_1+\phi_2) \right) - \partial_{\tau}  |\nabla_y (\phi_1+\phi_2)|^2 \right. & \notag \\ 
	\left. -(\gamma-1) \partial_{\tau}(\phi_1+\phi_2) \nabla_y^2 (\phi_1+\phi_2) +  2 c \partial^2_{\tau z} \psi -c^2 \nabla_y^2 \psi  \right] &=0, \label{eq: sys4} \\
	\varepsilon^2 \tilde{\varepsilon}^2 \left[- \partial_{\tau} \left( \partial_z (\phi_1+\phi_2) \right)^2 -(\gamma -1)  \partial_{\tau}(\phi_1+\phi_2) \partial_z^2 (\phi_1+\phi_2) - c^2 \partial_z^2 \psi \right] &=0. \label{eq: sys5}
\end{align}
The order of the system of equations should stay as they are stated above. Therefore, the relation $\varepsilon \leq \tilde{\varepsilon} \leq \sqrt{\varepsilon}$ should hold such that we arrive at different cases. The following table gives an overview on these. 
\renewcommand*{\arraystretch}{1.3}
	\begin{table}[h!]\label{tab:cases}
		\centering
	\begin{tabular}[h]{ll||c|c|c|c| c} \toprule
		Case & & equ. \eqref{eq: sys1} &  equ. \eqref{eq: sys2} &  equ. \eqref{eq: sys3} &  equ. \eqref{eq: sys4} &  equ. \eqref{eq: sys5} \\  
		\bottomrule \toprule
		Case 1: &  $\tilde{\varepsilon}=\varepsilon^{1/2}$& $\varepsilon^{3/2}$ &  $\varepsilon^{2}$ &  $\varepsilon^{2}$ &  $\varepsilon^{5/2}$ & $\varepsilon^{3}$  \\
		Case 2: &  $\varepsilon<\tilde{\varepsilon}<\varepsilon^{1/2}$ & 	$\varepsilon \tilde{\varepsilon}$ &  $\varepsilon^2 $ & $\varepsilon \tilde{\varepsilon}^2$ &  $\varepsilon^2 \tilde{\varepsilon}$ & $\varepsilon^2 \tilde{\varepsilon}^2$  \\
		Case 3:  &  $\tilde{\varepsilon}=\varepsilon$ & $\varepsilon^2$ &  $\varepsilon^2$ &  $\varepsilon^3$ &  $\varepsilon^3$ & $\varepsilon^4$  \\
		\bottomrule
	\end{tabular}
\vspace{0.3cm}
	\caption{Relation between $\varepsilon$ and $\tilde{\varepsilon}$ for  equations \eqref{eq: sys1} - \eqref{eq: sys5}}
\end{table}
\vspace{-0.3cm}

Note that the waves with velocity potentials $\phi_1$ and $\phi_2$ follow a parabolic approximation whereas the resulting field is omni-directional and the hydrophone is typically located far from the focal point, which means that the parabolic approximation is not ideal for the resulting wave that is here expressed in terms of $\psi$. So, since $\phi_1$,  $\phi_2$ are paraxial $\tilde{\varepsilon}$ is quite small from the point of view of $\phi_1$,  $\phi_2$  whereas for $\psi$ it holds that $\tilde{\varepsilon}$ is not that small from its point of view.

The following table gives an overview of the different models that we end up with. The subcases refer to a different accuracy or respectively to different orders of $\varepsilon$ up to which terms are considered. The enumeration (\lowroman{1}) and (\lowroman{2}) refers to the first and second equation in the following case distinctions, respectively. The colors emphasize that these cases are very similar and just differ in some constants.  
	\begin{table}[h!]
		\centering
	\begin{tabular}[h]{ll||c|c|c|c c|c c} \toprule
		Case & & equ. \eqref{eq: sys1} & equ. \eqref{eq: sys2} & equ. \eqref{eq: sys3}  & \multicolumn{2}{c}{equ. \eqref{eq: sys4} } & \multicolumn{2}{|c}{equ. \eqref{eq: sys5} }  \\ 
		& & & & & $\psi$ & $\phi_1+\phi_2$ & $\psi$ & $\phi_1+\phi_2$ \\ \bottomrule \toprule 
		\color{violet} Case 1: A &  $\tilde{\varepsilon}=\varepsilon^{1/2}$& (\lowroman{1}) & (\lowroman{2}) & (\lowroman{2}) & (\lowroman{2}) & &  (\lowroman{2}) & \\
		\color{cyan} Case 1: B &  $\tilde{\varepsilon}=\varepsilon^{1/2}$& (\lowroman{1}) & (\lowroman{2}) & (\lowroman{2}) & (\lowroman{2})  & (\lowroman{2}) &  (\lowroman{2}) & \\ \hline
		Case 2: A &  $\varepsilon<\tilde{\varepsilon}<\varepsilon^{1/2}$ & (\lowroman{1}) & (\lowroman{2}) & & (\lowroman{2}) &   & (\lowroman{2}) &  \\
		\color{violet} Case 2: B &  $\varepsilon<\tilde{\varepsilon}<\varepsilon^{1/2}$ & (\lowroman{1}) & (\lowroman{2}) & (\lowroman{2}) & (\lowroman{2}) &   & (\lowroman{2}) &  \\
		\color{cyan} Case 2: C &  $\varepsilon<\tilde{\varepsilon}<\varepsilon^{1/2}$ & (\lowroman{1}) & (\lowroman{2}) & (\lowroman{2}) & (\lowroman{2}) &  (\lowroman{2})  & (\lowroman{2}) &  \\ \hline
		Case 3:  &  $\tilde{\varepsilon}=\varepsilon$ & (\lowroman{1}) & (\lowroman{1}) & (\lowroman{2}) & (\lowroman{2}) & (\lowroman{2})  & (\lowroman{2}) & \\
		\bottomrule
	\end{tabular}
\vspace{0.3cm}
\caption{Overview of the following cases on the relation of $\varepsilon$ and $\tilde{\varepsilon}$ for  equations \eqref{eq: sys1} - \eqref{eq: sys5}}
	\end{table}

\underline{Case 1: A} For the case $\tilde{\varepsilon}=\sqrt{\varepsilon}$ considering only terms up to order $\varepsilon^2$ we end up with the model
\begin{align} 
	\frac{2 }{c}\partial^2_{\tau z} \phi_k - \nabla_y^2 \phi_k&=0, \qquad \text{ for } k \in \left\lbrace 1,2\right\rbrace \label{eq: Case1Ai} \\
	\tilde{\varepsilon} \left[ \frac{2}{c} \partial^2_{\tau z} \psi - \nabla_y^2 \psi\right] -\tilde{\varepsilon}^2 \partial_z^2 \psi  & =\frac{\gamma +1}{2c^4}\partial_{\tau}\left| \partial_{\tau} (\phi_1+\phi_2)\right|^2 + \partial_z^2 (\phi_1+\phi_2), \label{eq: Case1Aii}
\end{align}
where we split the first equation in $\phi_1 + \phi_2$ into two separate equations for $\phi_1$ and $\phi_2$ since they are linear and can be considered independently. Observe that the way equation \eqref{eq: Case1Aii} is written indicates the solution approach, in that equation \eqref{eq: Case1Ai} is solved for $\phi_1$, $\phi_2$, which then build a source term driving the second equation. \\

\underline{Case 1: B} For the case $\tilde{\varepsilon}=\sqrt{\varepsilon}$ considering terms up to order $\varepsilon^{5/2}$ we end up with the model
\begin{align} 
	\frac{2 }{c}\partial^2_{\tau z} \phi_k - \nabla_y^2 \phi_k&=0, \qquad \text{ for } k \in \left\lbrace 1,2\right\rbrace \label{eq: Case1Bi}\\
	\tilde{\varepsilon} \left[ \frac{2}{c} \partial^2_{\tau z} \psi - \nabla_y^2 \psi\right] -\tilde{\varepsilon}^2 \partial_z^2 \psi  & =\frac{\gamma +1}{2c^4}\partial_{\tau}\left| \partial_{\tau} (\phi_1+\phi_2)\right|^2 + \partial_z^2 (\phi_1+\phi_2) \notag \\
	& - \frac{\tilde{\varepsilon}}{c^2} \left[ \frac{2}{c} \left( \partial_{\tau}^2 (\phi_1+\phi_2)\partial_{z} (\phi_1+\phi_2)+ \gamma \partial_{\tau}(\phi_1+\phi_2) \partial_{\tau z}^2 (\phi_1+\phi_2) \right) \right. \notag \\
	&  \left. - \partial_{\tau}  |\nabla_y (\phi_1+\phi_2)|^2 -(\gamma-1) \partial_{\tau}(\phi_1+\phi_2) \nabla_y^2 (\phi_1+\phi_2) \right].
\end{align}
Using the first equation and plugging 
\begin{equation*}
	\nabla_y^2 \phi_k = \frac{2 }{c}\partial^2_{\tau z} \phi_k
\end{equation*}
into the second equation, it reduces to
\begin{align}
	\tilde{\varepsilon} \left[ \frac{2}{c} \partial^2_{\tau z} \psi - \nabla_y^2 \psi\right] -\tilde{\varepsilon}^2 \partial_z^2 \psi  & =\frac{\gamma +1}{2c^4}\partial_{\tau}\left| \partial_{\tau} (\phi_1+\phi_2)\right|^2 + \partial_z^2 (\phi_1+\phi_2) \notag \\
	& - \frac{\tilde{\varepsilon}}{c^2} \left[ \frac{2}{c} \left( \partial_{\tau}^2 (\phi_1+\phi_2)\partial_{z} (\phi_1+\phi_2)+  \partial_{\tau}(\phi_1+\phi_2) \partial_{\tau z}^2 (\phi_1+\phi_2) \right) \right. \notag \\
	&  \left. - \partial_{\tau}  |\nabla_y (\phi_1+\phi_2)|^2 \right]. \label{eq: Case1Bii}
\end{align}

\underline{Case 2: A} In case  $\varepsilon<\tilde{\varepsilon}<\varepsilon^{1/2}$ for terms up to order $\varepsilon^2$ the system of equations is
\begin{align} 
	\frac{2 }{c}\partial^2_{\tau z} \phi_k - \nabla_y^2 \phi_k&=0, \qquad \text{ for } k \in \left\lbrace 1,2\right\rbrace \label{eq: Case2Ai}\\
	\tilde{\varepsilon} \left[ \frac{2}{c} \partial^2_{\tau z} \psi - \nabla_y^2 \psi\right] -\tilde{\varepsilon}^2 \partial_z^2 \psi  & =\frac{\gamma +1}{2c^4}\partial_{\tau}\left| \partial_{\tau} (\phi_1+\phi_2)\right|^2. \label{eq: Case2Aii}
\end{align}
Note, that the model contains the KZK equation
\begin{equation*}
	\frac{2 }{c}\partial^2_{\tau z} \phi_k - \nabla_y^2 \phi_k = \frac{\gamma +1}{2c^4}\partial_{\tau}\left| \partial_{\tau} (\phi_1+\phi_2)\right|^2,
\end{equation*}
but in our case the quadratic nonlinearity is decoupled and appears as a source term for the equation in terms of $\psi$. \\

\underline{Case 2: B} In case  $\varepsilon<\tilde{\varepsilon}<\varepsilon^{1/2}$ for terms up to order $\varepsilon \tilde{\varepsilon}^2$ the system of equations is
\begin{align} 
	\frac{2}{c}\partial^2_{\tau z} \phi_k - \nabla_y^2 \phi_k&=0, \qquad \text{ for } k \in \left\lbrace 1,2\right\rbrace \label{eq: Case2Bi}\\
	\tilde{\varepsilon} \left[ \frac{2}{c} \partial^2_{\tau z} \psi - \nabla_y^2 \psi\right] -\tilde{\varepsilon}^2 \partial_z^2 \psi  & =\frac{\gamma +1}{2c^4}\partial_{\tau}\left| \partial_{\tau} (\phi_1+\phi_2)\right|^2 + \frac{\tilde{\varepsilon}^2}{\varepsilon} \partial_z^2 (\phi_1+\phi_2). \label{eq: Case2Bii}
\end{align}
Note that since $\varepsilon<\tilde{\varepsilon}<\varepsilon^{1/2}$ it holds that 
\begin{equation*}
	\varepsilon^3 < \varepsilon^2 \tilde{\varepsilon} < \varepsilon \tilde{\varepsilon}^2  < \varepsilon^2
\end{equation*}
and therefore, 
\begin{equation*}
	\varepsilon < \tilde{\varepsilon} <\frac{ \tilde{\varepsilon}^2}{\varepsilon }  < 1.
\end{equation*}

\underline{Case 2: C} In case  $\varepsilon<\tilde{\varepsilon}<\varepsilon^{1/2}$ for terms up to order $\varepsilon^2 \tilde{\varepsilon}$ the system of equations is
\begin{align} 
	\frac{2 }{c}\partial^2_{\tau z} \phi_k - \nabla_y^2 \phi_k&=0, \qquad \text{ for } k \in \left\lbrace 1,2\right\rbrace \label{eq: Case2Ci}\\
	\tilde{\varepsilon} \left[ \frac{2}{c} \partial^2_{\tau z} \psi - \nabla_y^2 \psi\right] -\tilde{\varepsilon}^2 \partial_z^2 \psi  & =\frac{\gamma +1}{2c^4}\partial_{\tau}\left| \partial_{\tau} (\phi_1+\phi_2)\right|^2 + \frac{\tilde{\varepsilon}^2}{\varepsilon} \partial_z^2 (\phi_1+\phi_2) \notag \\
	& - \frac{\tilde{\varepsilon}}{c^2} \left[ \frac{2}{c} \left( \partial_{\tau}^2 (\phi_1+\phi_2)\partial_{z} (\phi_1+\phi_2)+  \partial_{\tau}(\phi_1+\phi_2) \partial_{\tau z}^2 (\phi_1+\phi_2) \right) \right. \notag \\
	&  \left. - \partial_{\tau}  |\nabla_y (\phi_1+\phi_2)|^2 \right], \label{eq: Case2Cii}
\end{align}
where we reduced the second equation by plugging the first one into the second. \\

\underline{Case 3:}  For the case $\tilde{\varepsilon}=\varepsilon$ considering terms up to order $\varepsilon^3$ we end up with the model
\begin{align} 
	\frac{2 }{c}\partial^2_{\tau z}(\phi_1+\phi_2) - \nabla_y^2 (\phi_1+\phi_2)&=\frac{\gamma +1}{2c^4}\partial_{\tau}\left| \partial_{\tau} (\phi_1+\phi_2)\right|^2, \label{eq: Case3i} \\
	\tilde{\varepsilon} \left[ \frac{2}{c} \partial^2_{\tau z} \psi - \nabla_y^2 \psi\right] -\tilde{\varepsilon}^2 \partial_z^2 \psi  & = \partial_z^2 (\phi_1+\phi_2) \notag \\
& - \frac{\tilde{\varepsilon}}{c^2} \left[ \frac{2}{c} \left( \partial_{\tau}^2 (\phi_1+\phi_2)\partial_{z} (\phi_1+\phi_2)+ \partial_{\tau}(\phi_1+\phi_2) \partial_{\tau z}^2 (\phi_1+\phi_2) \right) \right. \notag \\
&  \left. - \partial_{\tau}  |\nabla_y (\phi_1+\phi_2)|^2 \right], \label{eq: Case3ii}
\end{align}
where we reduced the second equation by plugging the first one into the second. In this case the first equation in $\phi_1 + \phi_2$ can not be split since it is nonlinear in $\phi_1 + \phi_2$.

\subsection{Transformation into frequency domain}

The next step is to transfer the time domain formulations into frequency domain with the time harmonic ansatz 
\begin{align*}
	\phi_k(\tau, z, y) &= \hat{\phi}_k(z,y) e^{\imath \omega_k \tau}  \qquad \text{ for } k \in \left\lbrace 1,2\right\rbrace, \\
	\psi(\tau, z,y)&=\Re{ \hat{\psi}(z,y) e^{\imath(\omega_1-\omega_2)\tau}}.
\end{align*}
We will therefore look at the cases from above and deal with them separately. Note that we assume $\omega_1 > \omega_2 >0$. 
\\

\underline{Case 1: A} For equation \eqref{eq: Case1Ai} the time harmonic ansatz yields
\begin{align*}
	\frac{2}{c}\partial^2_{\tau z} \phi_k - \nabla_y^2 \phi_k = \frac{2}{c} \partial^2_{\tau z} \hat{\phi}_k e^{\imath \omega_k \tau}- \nabla_y^2 \hat{\phi}_k e^{\imath \omega_k \tau} =  \left( \frac{2  \imath \omega_k}{c} \partial_z \hat{\phi}_k - \nabla_y^2 \hat{\phi}_k\right) e^{\imath \omega_k \tau} = 0.
\end{align*}
Therefore, 
\begin{equation} \label{eq: eq1_fd}
	\frac{2 \imath \omega_k }{c} \partial_z \hat{\phi}_k-\nabla_y^2 \hat{\phi}_k=0 \qquad \text{ for } k \in \left\lbrace 1,2\right\rbrace.
\end{equation}
Now, going on to the second equation. Starting with the right hand side of equation \eqref{eq: Case1Aii} the time harmonic ansatz yields 
\begin{align*}
	\frac{\gamma +1}{2c^4} \partial_{\tau}\left| \partial_{\tau} (\phi_1+\phi_2)\right|^2 & +  \partial_z^2 (\phi_1+\phi_2) \\
	 = 	& \frac{\gamma+1}{2c^4} \partial_{\tau}  \left| \imath \omega_1 \hat{\phi}_1 e^{\imath \omega_1 \tau} + \imath \omega_2\hat{\phi}_2 e^{\imath \omega_2\tau} \right|^2  + e^{\imath \omega_1 \tau} \partial_z^2 \hat{\phi}_1  +  e^{\imath \omega_2\tau} \partial_z^2 \hat{\phi}_2  \\
	 =&\frac{\gamma+1}{2c^4}\left[ 2(\omega_1-\omega_2) \omega_1 \omega_2 \Re{ \imath \hat{\phi}_1 \overline{\hat{\phi}}_2 e^{\imath (\omega_1-\omega_2) \tau} } \right] + e^{\imath \omega_1 \tau} \partial_z^2 \hat{\phi}_1  +  e^{\imath \omega_2\tau} \partial_z^2 \hat{\phi}_2 .
\end{align*}
For the left hand side of equation \eqref{eq: Case1Aii} the time harmonic ansatz 
yields 
\begin{align*}
	\tilde{\varepsilon} \left[ \frac{2}{c} \partial^2_{\tau z} \psi - \nabla_y^2 \psi\right] -\tilde{\varepsilon}^2 \partial_z^2 \psi 
	= \Re{\left[ \tilde{\varepsilon} \left( \frac{2 \imath (\omega_1 - \omega_2)}{c} \partial_z \hat{\psi} -  \nabla_y^2 \hat{\psi} - \tilde{\varepsilon} \partial_z^2 \hat{\psi} \right) \right]  e^{\imath(\omega_1-\omega_2)\tau}}. 
\end{align*}
So, putting the second equation together yields
\begin{align*} 
	\Re \left\lbrace \left( \frac{2 \imath (\omega_1 - \omega_2)}{c} \partial_z \hat{\psi} -  \nabla_y^2 \hat{\psi} - \tilde{\varepsilon} \partial_z^2 \hat{\psi}  \right. \right.  - &  \left. \left.  \frac{\gamma+1}{2\tilde{\varepsilon} c^4}\left[ 2(\omega_1-\omega_2) \omega_1 \omega_2 \imath \hat{\phi}_1 \overline{\hat{\phi}}_2  \right] \right)e^{\imath (\omega_1-\omega_2) \tau} \right\rbrace \\
	&= \frac{1}{\tilde{\varepsilon}} \left( e^{\imath \omega_1 \tau} \partial_z^2 \hat{\phi}_1  +  e^{\imath \omega_2\tau} \partial_z^2 \hat{\phi}_2\right) .
\end{align*} \normalsize

\underline{Case 1: B} The first equation in frequency domain is already given in \eqref{eq: eq1_fd}. Analogously to the case above, the time harmonic ansatz for the second equation yields 
\small
\begin{align*} 
	&\Re \left\lbrace\left(  \frac{2 \imath (\omega_1 - \omega_2)}{c} \partial_z \hat{\psi}  - \nabla_y^2 \hat{\psi} - \tilde{\varepsilon} \partial_z^2 \hat{\psi}-(\omega_1-\omega_2) \left[ \frac{\gamma+1}{\tilde{\varepsilon} c^4} \omega_1 \omega_2 \imath \hat{\phi}_1 \overline{\hat{\phi}}_2 + \frac{2 }{c^2} \nabla_y \hat{\phi}_1 \nabla_y \overline{\hat{\phi}}_2 \right] \right)e^{\imath (\omega_1-\omega_2) \tau} \right\rbrace \\
	&= \frac{1}{\tilde{\varepsilon}} \left( e^{\imath \omega_1 \tau} \partial_z^2 \hat{\phi}_1  +  e^{\imath \omega_2\tau} \partial_z^2 \hat{\phi}_2 \right)  \\
	& - \frac{2}{c^3} \left[ -\omega_1^2 \partial_z \hat{\phi}_1^2 e^{2 \imath \omega_1 \tau} - \left( \left( \omega_1^2 + \omega_1 \omega_2\right)  \hat{\phi}_1 \partial_z \hat{\phi}_2 + \left( \omega_2^2 + \omega_1 \omega_2\right) \partial_z \hat{\phi}_1 \hat{\phi}_2\right) e^{\imath \left( \omega_1 + \omega_2\right) \tau} -\omega_2^2 \partial_z \hat{\phi}_2^2 e^{2 \imath \omega_2 \tau}  \right].
\end{align*} \normalsize

\underline{Case 2: A} Using the derivations in case 1, the model is here given by
\begin{align} 
	\frac{2 \imath \omega_k }{c} \partial_z \hat{\phi}_k-\nabla_y^2 \hat{\phi}_k&=0 \qquad \text{ for } k \in \left\lbrace 1,2\right\rbrace ,\\
	\frac{2 \imath (\omega_1 - \omega_2)}{c} \partial_z \hat{\psi} -  \nabla_y^2 \hat{\psi} - \tilde{\varepsilon} \partial_z^2 \hat{\psi}  &=\frac{\gamma+1}{\tilde{\varepsilon}c^4}\left[ (\omega_1-\omega_2) \omega_1 \omega_2 \imath \hat{\phi}_1 \overline{\hat{\phi}}_2  \right].
\end{align}
One can really observe here that the low frequency field propagates at the difference frequency $\omega_1-\omega_2$. 

\underline{Case 2: B} The system of equations only differs in a constant factor from case 1: A, so, we refer to the time harmonic model there. 

\underline{Case 2: C} The system of equations only differs in a constant factor from case 1: B. Therefore, we refer to case 1: B for the time harmonic model. \\

\underline{Case 3:} The time harmonic Ansatz yields for the first equation
\begin{align*}
	 \left( \frac{2 \imath \omega_1}{c} \partial_z \hat{\phi}_1 - \nabla_y^2 \hat{\phi}_1\right) e^{\imath \omega_1 \tau} +  \left( \frac{2  \imath \omega_2}{c} \partial_z \hat{\phi}_2- \nabla_y^2 \hat{\phi}_2\right) e^{\imath \omega_2 \tau} = \frac{\gamma+1}{c^4}\left[ (\omega_1-\omega_2) \omega_1 \omega_2 \Re{ \imath \hat{\phi}_1 \overline{\hat{\phi}}_2 e^{\imath (\omega_1-\omega_2) \tau}  } \right]
\end{align*}
and for the second
\begin{align*}
	&\Re \left\lbrace \left(  \frac{2\imath (\omega_1 - \omega_2)}{c} \partial_z \hat{\psi}  - \nabla_y^2 \hat{\psi} - \tilde{\varepsilon} \partial_z^2 \hat{\psi}- \frac{2(\omega_1-\omega_2)}{c^2}  \nabla_y \hat{\phi}_1 \nabla_y \overline{\hat{\phi}}_2 \right)e^{\imath (\omega_1-\omega_2) \tau} \right\rbrace \\
	&= \frac{1}{\tilde{\varepsilon}} \left( e^{\imath \omega_1 \tau} \partial_z^2 \hat{\phi}_1  +  e^{\imath \omega_2\tau} \partial_z^2 \hat{\phi}_2\right)  \\
	& - \frac{2}{c^3} \left[ -\omega_1^2 \partial_z \hat{\phi}_1^2 e^{2 \imath \omega_1 \tau} - \left( \left( \omega_1^2 + \omega_1 \omega_2\right)  \hat{\phi}_1 \partial_z \hat{\phi}_2 + \left( \omega_2^2 + \omega_1 \omega_2\right) \partial_z \hat{\phi}_1 \hat{\phi}_2\right) e^{\imath\left( \omega_1 + \omega_2\right) \tau} -\omega_2^2 \partial_z \hat{\phi}_2^2 e^{2\imath \omega_2 \tau}  \right].
\end{align*}

\subsection{Boundary conditions}

We consider a bounded computational domain $\tilde{\Omega}$ and due to the paraxial approach, where we assumed that the beams propagate in $z$-direction, we write $\tilde{\Omega}$ as a product
\begin{equation*}
	\tilde{\Omega}=\left(  0,L \right) \times \tilde{\Omega}_y
\end{equation*}
with $L>0$ and $\tilde{\Omega}_y$ a smooth subset of $\mathbb{R}^{d-1}$, so that $\left( z,y \right)  \in \left(  0,L\right)  \times \tilde{\Omega}_y$. 

Since the ultrasound beams are excited by a transducer, with $h_k(\tau, 0,y)=\hat{h}_k(y)e^{i \omega_k \tau}$ being the given time harmonic initial data at $z=0$ and absorbing boundary conditions on the boundary parts parallel to the $z$-axis we arrive at
\begin{align} \label{eq: bound5}
	\hat{\phi}_k(0,y)&=\hat{h}_k(y) \qquad \qquad \quad \text{ in } \tilde{\Omega}_y, \notag \\
	\partial_{\nu_y} \hat{\phi}_k(z,y)& = -\imath \sigma_k \hat{\phi}_k(z,y) \qquad (z,y) \in \left( 0, L \right) \times \partial \tilde{\Omega}_y.
\end{align}
Also for the velocity potential $\psi$ we assume impedance boundary conditions in order to avoid spurious reflections of waves on the boundary of the computational domain, i.e.
\begin{align}
	-\partial_z  \hat{\psi}(0,y) &= -\imath \sigma_0 \psih(0,y) \qquad \, \, \, y \in \tilde{\Omega} \notag \\
	\partial_z  \hat{\psi}(L,y) &= -\imath \sigma_L \psih(L,y) \qquad \, \, y \in \tilde{\Omega}, \\
	\partial_{\nu_y} \psih(z,y) &= -\imath \sigma \psih(z,y) \qquad \quad (z,y) \in \left( 0, L \right) \times \partial \tilde{\Omega}_y. \notag
\end{align}
Here $\sigma_0, \, \sigma_L >0$ and the coefficients $\sigma_k$ and $\sigma$ are nonnegative $L^{\infty}$ impedance coefficients that are bounded away from zero on an open subset of $\partial \tilde{\Omega}$. 

\subsection{Adaptions and generalizations of the model} \label{sec: adaptions}

\paragraph{Transducer setup.} There are several ways to arrange the two beams. In the considered model the transducer is split symmetrically in the middle, so along the central axis where the focus lies, such that the left and right side can be excited with different frequencies. Another possible arrangement could be to split the transducer symmetrically into three parts, so a broader part in the middle where one beam is excited and smaller parts on the left and right side for the other beams. It could also be interesting to use a random approach for the excitement of the transducer elements. 
In general it is as well possible to use two probes for the two different beams, but this involves further challenges on the application side which is the reason why medical tests are currently performed with only one probe.

\paragraph{Focusing.} A very important property of the beams in vibro-acoustography is their focusing. In the different models derived above we assume that the beams are naturally focused, so by geometrical means, which is assumed to be incorporated in the initial data $h_k$. Another way to achieve focusing is by adding time delays to achieve a sequence of delayed pulses which produce the equivalent of a lens, cf. e.g. \cite{szabo2013} . This is often referred to as electronic focusing. We briefly discuss the questions how electronic focusing can be accomplished and how it can be included in a paraxial model in vibro-acoustography. 

\begin{figure}[h!]
	\centering
	\includegraphics[width=0.5\textwidth]{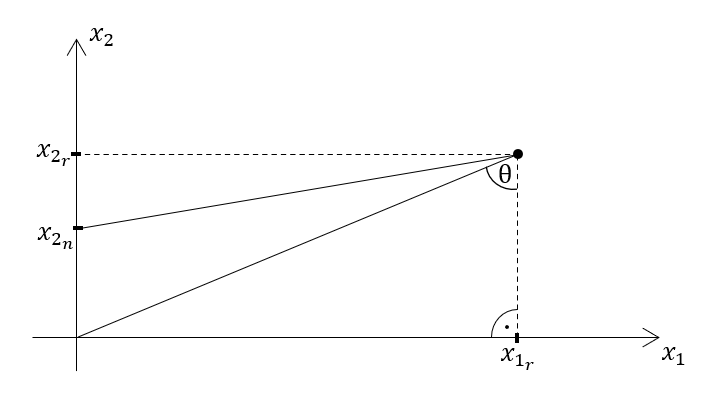}
	\caption{Electronic focusing of transducer elements}
	\label{fig: focus}
\end{figure}

Figure \ref{fig: focus} illustrates the geometry for the derivation of the time delay to focus an element of the transducer. In this case $x_1$ is the axis of propagation and the transducer is aligned to the $x_2$-direction. Furthermore, $x_{2_n}$ is the distance from the origin to the center of the $n$-th element of the transducer and at position $(x_{1_r}, x_{2_r})$ is the focus point.

The time delay to focus an element $x_{2_n}$ can be given by
\begin{equation*}
	\tilde{\tau}_n= \frac{1}{c} \left[ \sqrt{x_{2_r}^2 + x_{1_r}^2} - \sqrt{(x_{2_r} - x_{2_n})^2 + x_{1_r}^2}\right]. 
\end{equation*}
Under the assumptions that lateral variations are smaller than the axial distance, i.e. $x_{2_r} - x_{2_n} \ll x_{1_r}$ and the source is only excited at points $x_{2_r} \ll x_{1_r}$ the expression for $\tilde{\tau}_n$ can be approximated by
\begin{equation*}
	\tilde{\tau}_n \approx \frac{1}{c} \left[ \frac{x_{2_r}^2}{2x_{1_r}} - \frac{(x_{2_r} - x_{2_n})^2}{2x_{1_r}}\right]  = \frac{1}{c} \left[ -\frac{x_{2_n}^2}{2x_{1_r}} + x_{2_n} \frac{x_{2_r}}{x_{1_r}}\right] = \frac{1}{c} \left[ -\frac{x_{2_n}^2}{2x_{1_r}} + x_{2_n} \tan (\theta) \right].
\end{equation*}
Since $\cos(\theta)= \frac{x_{2_r}}{x_{1_r}} \approx 1$ and $\tan(\theta) \approx \sin(\theta)$, it follows that
\begin{equation} \label{eq:41}
	\tilde{\tau}_n \approx \frac{1}{c} \left[ -\frac{x_{2_n}^2}{2x_{1_r}} + x_{2_n} \sin (\theta) \right].
\end{equation}

\begin{figure}[h]
	\centering
	\includegraphics[width=0.5\textwidth]{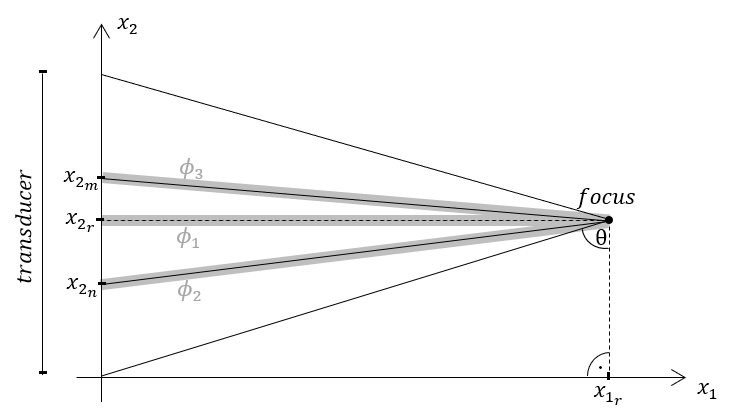}
	\caption{Focusing of three beams with time delays}
	\label{fig: focus3}
\end{figure}

For a focused beam at least three points are necessary. Figure \ref{fig: focus3} gives an overview of a model where we consider three velocity potentials $\phi_1$, $\phi_2$ and $\phi_3$. Here, $\phi_1$ starts at position $x_{2_r}$ and to keep it more general we assume that the focus point and $x_{2_r}$ are not on the central axis $x_1$ and we also do not assume symmetry. The other two beams start at position $x_{2_n}$ and $x_{2_m}$ respectively. 

The next step is to include these time delays. Since $\phi_1$ is at the center in our model no additional time delay is needed. Therefore, we perform the paraxial change of variables
\begin{equation}
	\phi_1(\tau, z, x) = \phi_1\left( t - \frac{x_{1}}{c}, \tilde{\varepsilon}x_1, \sqrt{\tilde{\varepsilon}}x_2\right) .
\end{equation}

The velocity potential $\phi_2$ starts at position $x_{2_n}$, which gives a time delay of $\tilde{\tau}_n$ and the direction of propagation is $\sqrt{(x_{2} - x_{2_n})^2 + x_{1}^2}$. Similar to what we have seen above when we derived the approximated time delay $\tilde{\tau}_n$, also the direction of propagation can be approximated, which yields
\begin{equation*}
	\sqrt{(x_{2} - x_{2_n})^2 + x_{1}^2} \approx x_{1} + \frac{(x_{2} - x_{2_n})^2}{2x_{1}}.
\end{equation*}
Therefore, we perform the following paraxial approach 
\begin{equation*} \label{eq: 43}
	\phi_2(\tau_n, z, x)=\phi_2 \left( t-\frac{x_{1}}{c}-\frac{(x_{2} - x_{2_n})^2}{2 c x_{1}} - \tilde{\tau}_n, \tilde{\varepsilon}x_1, \sqrt{\tilde{\varepsilon}}x_2\right) = \phi_2 \left( \tau- \tau_{r_n} - \tilde{\tau}_n, \tilde{\varepsilon}x_1, \sqrt{\tilde{\varepsilon}}x_2\right),
\end{equation*}
where $\tau_{r_n}=\frac{(x_{2} - x_{2_n})^2}{2 c x_{1}}$ and $\tilde{\tau}_n$ from \eqref{eq:41}. Similarly for the velocity potential $\phi_3$ we get
\begin{equation*} \label{eq: 44}
	\phi_3(\tau_m, z, x)=\phi_3 \left( t-\frac{x_{1}}{c}-\frac{(x_{2} - x_{2_m})^2}{2 c x_{1}} - \tilde{\tau}_m, \tilde{\varepsilon}x_1, \sqrt{\tilde{\varepsilon}}x_2\right) = \phi_3 \left( \tau- \tau_{r_m} - \tilde{\tau}_m, \tilde{\varepsilon}x_1, \sqrt{\tilde{\varepsilon}}x_2\right).
\end{equation*}
Now, analogously to above one could derive a model for vibro-acoustic imaging with electronic focusing.

\paragraph{Multiple beams.} Consider for example the model in Case 2: A given by equation \eqref{eq: Case2Ai}, \eqref{eq: Case2Aii}. This model can simply be extended to an arbitrary number of beams that are excited with two different frequencies $\omega_1, \, \omega_2$. Suppose that $n$ beams are excited with frequency $\omega_1$ and $m$ beams with $\omega_2$. The model in time domain is
\begin{align}
2 \partial^2_{\tau z} \phi_j -c \Delta_y\phi_j  &=0, \label{eq: 312} \\
\tilde{\varepsilon} \left[ \frac{2}{c} \partial^2_{\tau z} \psi - \nabla_y^2 \psi\right] -\tilde{\varepsilon}^2 \partial_z^2 \psi  &= \frac{\gamma +1}{2c^4} \partial_{\tau}\left| \partial_{\tau} \sum_{i=1}^{n+m} \phi_i \right|^2, \label{eq: 313}
\end{align}
for $j \in \left\lbrace 1, \ldots, n+m\right\rbrace$. In frequency domain with a time harmonic ansatz this yields
\begin{align}
\imath \omega_1 \partial_{ z} \hat{\phi}_k - c \Delta_y \hat{\phi}_k&= 0, \\
\imath \omega_2 \partial_{ z} \hat{\phi}_l - c \Delta_y \hat{\phi}_l&= 0, \\
\frac{2 \imath (\omega_1 - \omega_2)}{c} \partial_z \hat{\psi} -  \nabla_y^2 \hat{\psi} - \tilde{\varepsilon} \partial_z^2 \hat{\psi}  &= \frac{\gamma+1}{c^4} (\omega_1-\omega_2) \omega_1 \omega_2 \imath \left( \sum_{k=1}^n \hat{\phi}_k\right) \left( \sum_{l=n+1}^{n+m}\overline{ \hat{\phi}}_l\right) 
\end{align}
for $k \in \left\lbrace 1, \ldots n\right\rbrace$, $l \in \left\lbrace n+1, \ldots n+m\right\rbrace$ and arbitrary $m, n \in \mathbb{N}$. Taking multiple beams into account can be of interest to reduce the ill-posedness of the inverse problem.


\section{The inverse problem} \label{sec: IP}

The inverse problem of nonlinearity parameter imaging  will be considered for one of these models, namely Case 2A
\begin{equation}\label{PDEs_intro}
	\begin{aligned}
		&2\imath\omega_k\slo\partial_z\phih_k-\Delta_\zbot \phih_k = 0\mbox{ in }\tilde{\Omega} \quad  k\in\{1,2\}\\
		&2 \imath \omega_d \slo \partial_z\psih - \tilde{\varepsilon} \partial_z^2\psih - \Delta_\zbot \psih = 
		\imath\omega_p \tilde{\varepsilon}^{-1}\,\nlc \,\phih_1\, \overline{\phih}_2\mbox{ in }\tilde{\Omega},
	\end{aligned}
\end{equation}
on $\tilde{\Omega}=(0,L)\times\tilde{\Omega}_\zbot$, where $s=\frac{1}{c}$ and $\eta=\frac{\gamma+1}{c^4}$ with boundary conditions
\begin{equation}\label{BCs_intro}
	\begin{aligned}
		&\phih_k(0,\zbot)=\hat{h}_k(\zbot), \quad \zbot\in\tilde{\Omega}_\zbot\,,\\
		&\partial_{\nu_\zbot} \phih_k(z,\zbot)=
		-\imath\sigma_k \phih_k(z,\zbot) \quad (z,\zbot) \in(0,L)\times\partial\tilde{\Omega}_\zbot
		&&\quad k\in\{1,2\}\\
		&
		-\partial_z \psih(0,\zbot)=-\imath\sigma_0\psih(0,\zbot), \quad
		\partial_z \psih(L,\zbot)=-\imath\sigma_L\psih(L,\zbot), \quad \zbot\in \tilde{\Omega}_\zbot\\
		&\partial_{\nu_\zbot} \psih(z,\zbot)=
		-\imath\sigma \psih (z,\zbot) \quad (z,\zbot) \in(0,L)\times\partial\tilde{\Omega}_\zbot.
	\end{aligned}
\end{equation}

The parameter we want to identify is $\eta$, which is directly connected to the space dependent nonlinearity parameter $B/A$, since
\begin{equation*}
	\eta = \frac{\gamma+1}{c^4} = \frac{2 \beta}{c^4} = \frac{2+ B/A}{c^4}.
\end{equation*}
In well-known nonlinear wave equations the nonlinearity parameter is often denoted by $\beta$ in the context of liquids or $\gamma$ for gases. 

We assume the speed of sound $c$ to be constant. Identifying the nonlinearity parameter and the speed of sound simultaneously leads to a nonlinear inverse problem that is dealt with in \cite{2023simultaneous}.



Since the excited low-frequency field in terms of $\psi$ propagates omni-directionally, alternatively to the paraxial form of the partial differential equation in terms of $\psih$ one can also consider the Helmholtz equation

	\begin{equation}\label{PDEsTransf}
		\begin{aligned}
			&2\imath\omega_k\slo\partial_z\phih_k-\Delta_\zbot \phih_k = 0\mbox{ in }\tilde{\Omega} \quad  k\in\{1,2\}\\
			&-\omega_d^2 \slo^2 \check{\psi} - \Delta_x \check{\psi} = \imath\omega_p \, 1_{\Omega}\Bigl(\check{\nlc}\, P^{-1}\bigl(\phih_1\, \overline{\phih}_2\bigr)\Bigr)\mbox{ in }\Omega
		\end{aligned}
	\end{equation}
	for $\check{\psi}=P^{-1}\psih$, $\check{\nlc}=P^{-1}\nlc$ (cf. \eqref{parax_deriv}) with boundary conditions
	\begin{equation}\label{BCsTransf}
		\begin{aligned}
			&\phih_k(0,\cdot)=\hat{h}_k \mbox{ in }\tilde{\Omega}_\zbot\,,
			\quad \partial_{\nu_\zbot} \phih_k(z,\zbot)=
			-\imath\sigma_k \phih_k & \mbox{ in }(0,L)\times\partial\tilde{\Omega}_\zbot
			\quad k\in\{1,2\}\\
			&\partial_\nu \check{\psi}=-\imath\sigma \check{\psi}  \mbox{ on }\partial\Omega 
		\end{aligned}
	\end{equation}
	with $P(\Omega)\supseteq\tilde{\Omega}=(0,L)\times\tilde{\Omega}_\zbot$, where $P$ is the paraxial transform defined by \eqref{parax} and 
	$1_{\Omega}$ the extension by zero operator for a function defined on $P^{-1}(\tilde{\Omega})$ to all of $\Omega$. This would allow to use a possibly larger propagation domain for the $\psi$ wave as compared to the one for the $\phi_1$, $\phi_2$ waves, but come at the price of a more involved implementation. 

In the following we deal with the inverse problem of nonlinearity imaging with the alternative formulation in \eqref{PDEsTransf}. Since $\check{\eta}$ only appears in the second equation, only this equation is relevant for reconstruction, so that the system of equations reduces to an inhomogeneous Helmholtz equation and $\phih_1, \, \phih_2$ can be computed independently beforehand. \\  

\subsection{Formulation of the inverse problem}

For the case of $\slo$ being known and constant, consider the reconstruction of $\check{\nlc}=\check{\nlc}(x)$ in 
\begin{equation}\label{bvp_psicheck}
	\begin{aligned}
		-\omega_{d}^2 \slo^2 \check{\psi} - \Delta_x \check{\psi} = \check{\nlc}\, f \qquad &  \mbox{ in }\Omega\\
		\partial_\nu \check{\psi}=-\imath\sigma \check{\psi} \qquad & \mbox{ on }\partial\Omega 
	\end{aligned}
\end{equation}
with $f=\imath\omega_{p} \, P^{-1}\bigl(\phih_{1}\, \overline{\phih}_{2}\bigr)$  
cf. \eqref{PDEsTransf} from observations

\begin{equation} \label{observations}
	\check{p}^{\textup{meas}}=\imath \omega_{d}\mbox{tr}_{\check{\Gamma}}\check{\psi}\,,
\end{equation}

where $\check{\Gamma}$ is a manifold representing the receiver array immersed in the acoustic domain $\Omega$. 

This inverse problem can be written as a linear operator equation
\begin{equation} \label{lip_op}
	T\check{\nlc} = y
\end{equation}

with $y= \check{p}^{\textup{meas}}$ and with the forward operator $T=C \circ S$ being a concatenation of the parameter-to-state map
\begin{equation*}
	S: \check{\nlc} \mapsto \check{\psi} \text{ solving \eqref{bvp_psicheck}}
\end{equation*}
with the (linear) observation operator 
\begin{equation*}
	C: \check{\psi} \mapsto \imath \omega_d \mbox{tr}_{\check{\Gamma}} \check{\psi}.
\end{equation*}
Note that since $(\phih_1,\phih_2)$ can be computed independently of $\eta$, the forward operator $S$ is linear. So, this means that the linear operator $T$ is given by 
\begin{equation} \label{forwOp}
	T:L^2(\Omega)\to L^2(\check{\Gamma}), \quad 
	\check{\nlc}\mapsto \imath \omega_{d}\mbox{tr}_{\check{\Gamma}}\check{\psi}, 
\end{equation}
 where $\check{\psi}$ solves \eqref{bvp_psicheck}.
Well-definedness and continuity of the linear operator $T$ for $f\in L^p(\Omega)$ 
follows from known results on well-posedness of the Helmholtz equation with impedance boundary conditions; cf. e.g. \cite{ThesisMelenk, MoiolaSpence:2014} and the references therein.  

\subsection{Solving the inverse problem}

For solving the inverse problem we want to use the Landweber-Kaczmarz method. A Landweber step for solving 
\eqref{lip_op} is defined by a gradient descent step for the least squares functional
\begin{equation*}
	\left\| T(\check{\nlc}^{(n)})-y\right\|^2_{L^2}
\end{equation*}
i.e., by 
\begin{equation} \label{Landweber}
	\check{\nlc}^{(n+1)}= \check{\nlc}^{(n)}+\mu T^*(y-T(\check{\nlc}^{(n)})).
\end{equation}
with an appropriately chosen step size $\mu$. Therefore, we want to derive the adjoint operator of $T(\check{\nlc})$.  

Observe that for all $\check{\eta} \in L^2(\Omega), \, q \in L^2(\check{\Gamma})$ it holds that
\begin{equation*}
	\left\langle T\check{\nlc},q \right\rangle_{L^2(\check{\Gamma})} = \int_{\check{\Gamma}} T\check{\nlc} \cdot \overline{q} \dd{\check{\Gamma}} = \int_{\check{\Gamma}} \imath \omega_{d} \check{\psi} \cdot \overline{q} \dd{\check{\Gamma}}= \int_{\Omega} \check{\nlc} \cdot T^*q \dd{x} = \left\langle \check{\nlc}, T^* q\right\rangle_{L^2(\Omega)}
\end{equation*}
and where $\check{\psi} \in H^1(\Omega; \mathbb{C})$ satisfies the weak direct equation 
\small
\begin{align*}
	 \int_{\partial \Omega} \imath \sigma \check{\psi} \overline{w} \dd{S} + \int_{\Omega} \nabla \check{\psi} \nabla \overline{w} \dd{x} - \int_{\Omega}\omega_{d}^2 s^2 \check{\psi} \overline{w} \dd{x} = \int_{\Omega} \check{\nlc} \, f\overline{w} \dd{x} \qquad \text{for all } w \in H^1(\Omega; \mathbb{C}).
\end{align*}
\normalsize 
Thus $T^*q=\overline{f}r$ where $r$ solves the weak form of the adjoint equation
\begin{equation*}
	\int_{\partial \Omega \setminus \check{\Gamma}} \imath \sigma v \overline{r} \dd{S} +\int_{\Omega} \nabla v \nabla \overline{r} \dd{x} - \int_{\Omega} \omega_{d}^2 s^2 v \overline{r} \dd{x} = \int_{\check{\Gamma}} \imath \omega_{d}v \overline{q} \dd{\check{\Gamma}} \qquad \text{for all } v \in H^1(\Omega; \mathbb{C})
\end{equation*}
and the adjoint operator is given by
\begin{equation} \label{adjOp}
	T^*: L^2 (\check{\Gamma}) \to L^2(\Omega), \quad q  \mapsto f \overline{r},
\end{equation}
where $r$ solves
\begin{align} \label{PDEadj}
	-\omega_d^2 s^2 r - \Delta r &= 0 \quad \mbox{ in } \Omega \notag \\
	\partial_\nu r &= \imath \sigma r \quad   \mbox{ on } \partial \Omega \setminus \check{\Gamma} \\
	\left[ \partial_\nu r \right]_{\check{\Gamma}} = -\imath \omega_d q \quad \text{ on } \check{\Gamma} \cap \Omega, & \qquad \partial_\nu r = -\imath \omega_d q \quad \text{ on } \check{\Gamma} \cap \partial \Omega, \notag
\end{align}
where $\left[ \partial_\nu r \right]$ denotes the jump of the normal derivative over the interface $\check{\Gamma}\cap \Omega$ (in case $\check{\Gamma}$ or part of it is immersed in the computational domain $\Omega$). 

Finally, to recover $\check{\nlc}$ we take into account measurements at multiple frequencies $\omega_{1,j}$, $\omega_{2,j}$, 
$j\in  \left\lbrace 1, \ldots, J\right\rbrace$ and also several receiver array locations $\check{\Gamma}^m$, $m \in \left\lbrace 1, \ldots, M\right\rbrace$ might be used.
Note that for elementary dimensionality reasons a general function $\check{\nlc}=\check{\nlc}(x)$ cannot be expected to be uniquely determined from a single boundary observation. 

Thus, we actually deal with a set of several model and observation operators, which gives the forward operator $T_p=C^m \circ S^j$ and data $y_p$ for $p=(m-1)J+j$. So, we can write the inverse problem of reconstructing $\check{\nlc}$ as a system of operator equations
\begin{equation*}
	T_p(\check{\nlc}) =y_p, \qquad p \in \left\lbrace 1, \ldots, P=J\cdot M \right\rbrace 
\end{equation*}
and apply Kaczmarz type method by sequentially performing one Landweber step in a cyclically repeated manner
\begin{equation} \label{Landweber-Kaczmarz}
	\check{\nlc}^{(n+1)}= \check{\nlc}^{(n)}+\mu T_p^*(y_p-T_p(\check{\nlc}^{(n)})),
\end{equation}
where $p=\mod(n-1,P)+1$.


\section{Algorithms} \label{sec: num}

We present algorithms for solving the direct and also the inverse problem for the model in Case 2A given by \eqref{eq: Case2Ai}, \eqref{eq: Case2Aii}. Therefore, the following modules are used as subroutines to the algorithms: 

\begin{itemize}
	\item $\phih_k=$solveParaxialPDE$(\omega_k, s, \sigma_k, h_k)$: solve
	\begin{align} \label{paraxialPDE}
		2\imath\omega_k\slo\partial_z\phih_k-\Delta_\zbot \phih_k= 0 \qquad &\mbox{ in }\tilde{\Omega}, \notag \\
		\phih_k(0,\cdot)=\hat{h}_k \qquad & \mbox{ in }\tilde{\Omega}_\zbot\,, \\
		\partial_{\nu_\zbot} \phih_k(z,\zbot)=
		-\imath\sigma_k \phih_k \qquad & \mbox{ in }(0,L)\times\partial\tilde{\Omega}_\zbot;  \notag
	\end{align}
	\item $\check{\psi}=$solveHelmholtz$(\omega_k, s, \sigma, \check{\eta},f)$: solve
	\begin{align*}
		-\omega_{d}^2 \slo^2 \check{\psi} - \Delta \check{\psi} = \check{\nlc}\, f \qquad &  \mbox{ in }\Omega, \\
		\partial_\nu \check{\psi}=-\imath\sigma \check{\psi} \qquad & \mbox{ on }\partial\Omega,
	\end{align*}
 with $f=\imath\omega_{p} \, P^{-1}\bigl(\phih_{1}\, \overline{\phih}_{2}\bigr)$;
	\item $r=$solveHelmholtz$\Gamma(\omega_k, s, \sigma, q)$: solve
	\begin{align*} 
	-\omega_d^2 s^2 r - \Delta r &= 0 \quad \mbox{ in } \Omega \notag \\
	\partial_\nu r &= \imath \sigma r \quad   \mbox{ on } \partial \Omega \setminus \check{\Gamma} \\
	\left[ \partial_\nu r \right]_{\check{\Gamma}} = -\imath \omega_d q \quad \text{ on } \check{\Gamma} \cap \Omega, & \qquad \partial_\nu r = -\imath \omega_d q \quad \text{ on } \check{\Gamma} \cap \partial \Omega, \notag
	\end{align*}
\end{itemize}
The submodules $\check{\psi}=$solveHelmholtz$(\omega_k, s, \sigma, \check{\eta},f)$ and $r=$solveHelmholtz$\Gamma(\omega_k, s, \sigma, q)$ can be implemented using standard solvers for Helmholtz equation.

For the submodule $\phih_k=$solveParaxialPDE$(...)$ an individual algorithm has been designed. We assume that the beams are excited by a linear transducer, which means that in one dimension the transducer can be represented as $\tilde{\Omega}_k=[-l_k, l_k] \subseteq \tilde{\Omega}_y$ with $l_k >0$. We can make use of this geometry and consider the corresponding eigenvalue problem for the Laplace operator $\Delta_y$ with Neumann boundary conditions, i.e.
\begin{align*}
	-\Delta_{y} v_k^j(y)=\lambda_k^j v_k^j (y) &\qquad y \in \tilde{\Omega}_k, \\
	\partial_{\nu} v_k^j(y)=0 & \qquad y \in \partial  \tilde{\Omega}_k.
\end{align*}
for $k=1,2$, $j=1,2,\ldots$. Then, the explicit solution is given by
\begin{equation*}
	v_k^j(y)=\sqrt{\frac{1}{l_k}} \cos(j \pi \frac{y}{2 \,l_k}), \qquad \lambda_k^j=\frac{\pi^2 j^2}{4\,l_k^2}.
\end{equation*}
Generalization to higher dimensions can be achieved easily, since on a rectangular domain separation of variables is possible, e.g. for two dimensions
\begin{equation*}
	v^{mn}(y_1,y_2)=v^m(y_1) v^n(y_2), \qquad \lambda^{mn}=\lambda^m+\lambda^n
\end{equation*}
where $v^m(y_1)$ and $\lambda^m$ are the eigenfunctions and eigenvalues corresponding to the one-dimensional problem on $[-b_1,b_1 \,]$, respectively,  $v^n(y_2)$ and $\lambda^n$ are those corresponding to the one-dimensional problem on$[-b_2,b_2 \,]$ with $b_1, \, b_2 \, >0$. Furthermore, this can simply be adapted to other boundary conditions and geometries such as disks or spheres, where the explicit eigensystem of the eigenvalue problem of the Laplace operator is known, cf. \cite{Grebenkov_2013}. 
 
Overall, the velocity potential $\hat{\phi}_k$ at any position $(z,y)$ is given by
\begin{equation} \label{eq: 405}
	\hat{\phi}_k(z,y)=\sum_j a_k^j(z)v_k^j(y).
\end{equation}
Now, plugging this Ansatz into equation \eqref{paraxialPDE} yields
\begin{equation*}
	\sum_j \left( 2 \imath \omega_k s a_z^j(z) + \lambda^j a^j(z)\right) v^j_k(y)=0,
\end{equation*}
which means that we have to solve the ODE
\begin{equation*}
	 a_z^j=\imath \frac{\lambda_k^j}{ 2\omega_k s} a^j.
\end{equation*}
The explicit solution is then given by
\begin{equation}
	a^j = a^j(0) e^{\imath \frac{\lambda_k^j}{2 \omega_k s}z}.
\end{equation}
with
\begin{equation*}
	a^j(0)=\frac{1}{2l_k} \int_{-l_k}^{l_k} \hat{\phi}_k(0,y) v^j_k(y) \dd{y}.
\end{equation*}
Therefore, $\phih_k$ is given by
\begin{align*}
	\hat{\phi}_k(z,y)=
	\sum_j\frac{1}{2 l_k} \int_{-l_k}^{l_k} \hat{\phi}_k(0,y) v^j_k(y) \dd{y} e^{\imath \frac{\lambda_k^j}{2 \omega_k s}z} v_k^j(y), \qquad k \in \left\lbrace 1,2\right\rbrace . 
\end{align*}
In the discrete case calculating the initial data corresponds to applying a discrete cosine transform (DCT) and using the Ansatz means applying the inverse discrete cosine transform (IDCT). 


	\begin{table}[h!]\label{tab:cases}
		\centering
		\begin{tabular}[h]{l l} \toprule
			\textbf{Input:} $\hat{h}_k(0,y)$ & \% transducer excitation \\
			\textbf{Output:} $\phih_k(z_\text{out},y)$ & \\ 
			\hline
			$a^j_k(0) = \text{DCT}(\hat{h}_k(0, y))$  \qquad \qquad &	\% Ansatz for linear transducer \\
			for $z=0:\Delta z: z_\text{out}$ & \\
			\qquad $a^j_k(z) = a^j_k(0) e^{\imath \frac{\lambda_k^j}{2 \omega_k s}z}$ & \\
			$\phih_k ( z_\text{out}, y)=\text{IDCT}(a^j_k( z_\text{out}))$ \qquad \qquad & \\
		\bottomrule
		\end{tabular}
		\vspace{0.3cm}
		\caption{Pseudo-Code for $\phih_k=$solveParaxialPDE$(\omega_k, s, \sigma_k, \hat{h}_k)$}
	\end{table}

\subsection{Direct problem}

The direct problem consists of solving \eqref{PDEsTransf}, \eqref{BCsTransf} for $\check{\psi}$, where $\check{\eta}$ is known. The initial transducer signal is given by $\hat{h}_k$ and designed such that the beam focuses, which can be achieved by time delayed pulses, see section \ref{sec: adaptions} for more details. 

\begin{table}[h!]\label{tab:cases}
	\centering
	\begin{tabular}[h]{l l} \toprule
		\textbf{Input:} $\hat{h}_k(0, y)$ & \% transducer excitation \\
		\textbf{Output:} $\check{\psi}(x)$ & \\ 
		\hline
		for $z=0:\Delta z: z_{\text{out}}$ & \\
		\qquad $\phih_k(z, \tau)=$solveParaxialPDE$(\omega_k, s, \sigma_k, h_k)$ \qquad & \% paraxial variables \\
		$\check{\phi}_{12}=P^{-1}\bigl(\phih_{1}\, \overline{\phih}_{2}\bigr)$ & \% change of variables \\
		for $x=x_{\text{in}}:\Delta x: x_{\text{out}}$ & \\
		\qquad $\check{\psi}(x)=$solveHelmholtz$(\omega_k, s, \sigma, \check{\eta},\check{\phi}_{12})$  &	\% original variables \\
		
		\bottomrule
	\end{tabular}
	\vspace{0.3cm}
	\caption{Pseudo-Code for the forward problem}
\end{table}

\subsection{Inverse problem}

For the inverse problem we want to reconstruct $\check{\eta}$ in \eqref{bvp_psicheck} from given observations, cf. \eqref{observations}. Since $\check{\eta}$ only appears in the equation for $\check{\psi}$, we can compute $\phih_2, \, \phih_2$ separately as presented above, apply the inverse operator for the paraxial change of variables and use this as an input for the inhomogeneous Helmholtz equation in the following system.  

\begin{table}[h!]\label{tab:cases}
	\centering
	\begin{tabular}[h]{l l} \toprule
		\textbf{Input:} $\check{\eta}^{(n)}, \, \check{\phi}_{12}$ &  \\
		\textbf{Output:} $\check{\eta}^{(n+1)}$ & \\ 
		\hline
		\qquad $\check{\psi}^{(n)}=$solveHelmholtz$(\omega_k, s, \sigma, \check{\eta},\check{\phi}_{12})$  \qquad &	\% computation of $\check{\psi}=S(\check{\eta}^{(n)})$ \\
		\qquad $q=y-\mbox{tr}_{\check{\Gamma}}\check{\psi^{(n)}}$ & \% computation of adjoint states \\
		\qquad $r^{(n)}=$solveHelmholtz$\Gamma(\omega_k, s, \sigma, q)$ & \\
		\qquad $\xi^{(n)}=\imath \omega_p \check{\phi}_{12} r^{(n)} $ & \\
		
		set $\check{\eta}^{(n+1)}=\check{\eta}^{(n)} + \mu \xi^{(n)}$ & \\
		
		\bottomrule
	\end{tabular}
	\vspace{0.3cm}
	\caption{Pseudo-Code for the inverse problem}
\end{table}

\section{Outlook}

In this paper we have introduced a new approach for directive models in vibro-acoustic imaging, where we used a paraxial change of variables. We have transferred the systems of equations to frequency domain, studied the identification of the nonlinearity parameter for one of these models, and derived a Landweber-Kaczmarz method for reconstruction. 

Concerning forward simulation we have presented a numerical method that takes the geometry of the transducer into account and explicitly solves the paraxial partial differential equation. Next steps are to implement the devised algorithms for the direct problem and the Landweber-Kaczmarz method for the reconstruction of the nonlinearity parameter. 

Future work would also consist of adapting the numerical method for the direct problem and also the inverse problem for variable speed of sound and variable nonlinearity parameter. While the reconstruction of the nonlinearity parameter alone yields a linear inverse problem where the first equation of the system in $\phi_1$, $\phi_2$ can be computed independently and the inverse problem then simply consists of the second equation, for the simultaneous reconstruction of both parameters the whole system has to be taken into account, which results in a nonlinear inverse problem.

\section*{Acknowledgment} The author wishes to thank her supervisor Barbara Kaltenbacher, University of Klagenfurt, for fruitful and inspiring discussions. This work was supported by the Austrian Science Fund (FWF) under the grant DOC78.








\end{document}